\documentclass[12pt,eqno]{article}
\usepackage{amsmath,amsthm,amscd,amssymb}
\usepackage{latexsym}
\newtheorem{remark}{Remark}[section]
\newcommand{\bremark}{\begin{remark} \em}
\newcommand{\eremark}{\end{remark} }

\newcommand{\be}{\beta}
\newcommand{\la}{\lambda}
\newcommand{\om}{\omega}

\newcommand{\Om}{\Omega}

\begin{document}
\parindent 15pt
\renewcommand{\theequation}{\thesection.\arabic{equation}}
\renewcommand{\baselinestretch}{1.15}
\renewcommand{\arraystretch}{1.1}
\def\disp{\displaystyle}
\title{\large PPW and Chiti type inequalities for the eigenvalue problem of Robin Laplacian
\footnotetext{This work is supported by  NNSFC(No.10971061 and No.11271120) and by Hunan Provincial Innovation
Foundation For Postgraduate (No.CX2011B198)\\
E-mail:qiuyidai@yahoo.cn,  shifeilin1116@163.com} \footnotetext{}}
\author{{\small Qiuyi Dai, \hspace{1cm}Feilin Shi }\\
{\small Department of Mathematics, Hunan Normal University}\\
 {\small Changsha Hunan 410081, P.R.China}\\
}
\maketitle

\abstract{\small Let $\Omega \subset R^n(n\geq 2)$ be a bounded
domain with boundary $\partial\Om$, $\nu$ be the outward unit vector normal to $\partial\Om$, and $0<\beta< +\infty$ be a parameter. We prove two results for the following Robin eigenvalue problem
\begin{eqnarray*}\label{eq1.1}
\left\{
\begin{array}{ll}
-\Delta\psi=\lambda\psi &x\in\Omega,\\
\frac{\partial\psi}{\partial\nu}+\be\psi=0 &x\in\partial\Omega.
\end{array}
\right.
\end{eqnarray*}
One is an upper bound for the ratio of the first two eigenvalues which can be used to recover the PPW conjecture proved by M.S.Ashbaugh and R.D.Benguria in \cite{Ben} and  \cite{Ashb}, the other is a reverse H$\ddot{o}$lder inequality for the first eigenfunction which is a natural generalization of Chiti's reverse H$\ddot{o}$lder inequality for the first eigenfunction of Dirichlet Laplacian.
}
\vskip 0.1cm
{\bf AMS subject classification:} 35P15, 35P30, 35J65, 35J70
\vskip 0.1cm
{\bf Key words:} Robin Laplacian, decreasing rearrangement, eigenvalue problem

\section*{1.\ Introduction}

\setcounter{section}{1}

\setcounter{equation}{0}

\noindent Let $\Omega \subset R^n(n\geq 2)$ be an open bounded
domain whose boundary $\partial\Om$ is assumed to be of Lipschitz
type. We consider the following eigenvalue problem
\begin{eqnarray}\label{eq1.1}
\left\{
\begin{array}{ll}
-\Delta\psi=\lambda\psi &x\in\Omega,\\
\frac{\partial\psi}{\partial\nu}+\be\psi=0 &x\in\partial\Omega,
\end{array}
\right.
\end{eqnarray}
where $\Delta=\sum_{i=1}^{n}\frac{\partial^2}{\partial x_{i}^2}$ is the Laplace operator,and $0\leq\beta\leq+\infty$ is a parameter.

It is well known that problem (\ref{eq1.1}) has a purely discrete real spectrum $\{\lambda_k(\Omega,\beta)\}_{k=1}^{+\infty}$ which can be arranged in an increasing way as the following
$$0\leq\lambda_1(\Om,\be)<\lambda_2(\Om,\be)\leq\lambda_3(\Om,\be)\leq\cdots\leq\lambda_k(\Om,\be)\rightarrow+\infty,\ k\rightarrow+\infty.$$
Here each eigenvalue is repeated according to its multiplicity.

The study of eigenvalue problems has its fundamental importance in mathematical physics and mathematics itself. Much attention has been paid to the estimate of the eigenvalues, as well as of the norm of eigenfunctions, and many results have been derived for the special cases $\be=0$ and $\be=+\infty$ of problem (\ref{eq1.1}) (see for example\cite{Ben,Ashb,Ashbaugh eigenvalue ratios,Ashbaugh-Neumann,Bandle,ref GC,G. Chiti2, Fab,Henrot1,Henrot2,ref SK,Kra,Kroger,PayRay,PS,Jobin,Weinberger}). We will mention some of these results which are closely related to our purpose of this paper in the following paragraphs.

When $\be=0$, problem (\ref{eq1.1}) is reduced to the following
\begin{eqnarray}\label{eq12.1}
\left\{
\begin{array}{ll}
-\Delta\psi=\lambda\psi &x\in\Omega,\\
\frac{\partial\psi}{\partial\nu}=0 &x\in\partial\Omega,
\end{array}
\right.
\end{eqnarray}
which is called Neumann eigenvalue problem for Laplace operator, or eigenvalue problem for Neumann Laplacian. It is easy to see that $\lambda_1(\Om,0)=0$ and the first nonzero eigenvalue of problem (\ref{eq12.1}) is $\lambda_2(\Om,0)$. For the simplicity of the notation, we traditionally denote $\lambda_k(\Om,0)$ by $\mu_{k-1}(\Om)$ for any $k\geq 1$. Let $\Om^*$ be the Schwarz symmetrization of $\Om$, that is, $\Om^*$ be the ball in $R^n$ with center at origin and such that $\Om^*$ and $\Om$ have the same volume. The most beautiful and important result is the following Szeg$\ddot{o}$-Weinberger inequality
\begin{eqnarray}\label{eq1.30.1}
\mu_1(\Om)\leq\mu_1(\Om^*)\ \ \ \mbox{with equality if and only if}\ \ \Omega\ \ \mbox{is a ball,}
\end{eqnarray}
which was proved by Szeg$\ddot{o}$ for dimension $n=2$ in \cite{szego}, and by Weinberger for dimension $n>2$ in \cite{Weinberger}. Some more results about problem (\ref{eq12.1}) can be found in \cite{Ashbaugh-Neumann,Kroger} etc. We also remark here that
\begin{eqnarray}\label{eq12.2}
\lambda_1(\Om,\be)\rightarrow 0\ \ \ \mbox{and}\ \ \ \lambda_2(\Om,\be)\rightarrow\mu_1(\Om)>0\ \ \mbox{as}\ \ \beta\rightarrow 0^{+}
\end{eqnarray}
for any $\Omega\subset R^n$. Hence, (\ref{eq1.30.1}) implies that there exist a constant $\beta_0>0$ which maybe depends on $\Om$ such that
\begin{eqnarray}\label{eq12.3}
\lambda_2(\Om,\beta)\leq\lambda_2(\Om^*,\beta)
\end{eqnarray}
for any $0<\beta\leq\beta_0$ provided that $\Om$ is not a ball.

When $\be=+\infty$, problem (\ref{eq1.1}) is reduced to the following eigenvalue problem
\begin{eqnarray}\label{eq12.4}
\left\{
\begin{array}{ll}
-\Delta\varphi=\lambda\varphi &x\in\Omega,\\
\varphi=0 &x\in\partial\Omega,
\end{array}
\right.
\end{eqnarray}
which is called eigenvalue problem for Dirichlet Laplacian. As usual, we denote by $\lambda_k(\Om)$ the $k^{th}$ eigenvalue of the problem (\ref{eq12.4}).  Problem (\ref{eq12.4}) was extensively studied by many authors, and many interesting and important results were obtained (see \cite{Ben,Ashb,Ashbaugh eigenvalue ratios,Bandle,ref GC, G. Chiti2,Fab,Henrot1,Henrot2,ref SK,Kra,PayRay,PS,Jobin}). It is impossible to exhaust all results about problem (\ref{eq12.4}) in a small paper. Here, we restate some of them to motivate our purpose of the present paper. The first result we recall here is the following Faber-Krahn inequality
\begin{eqnarray}\label{eq12.5}
\lambda_1(\Om)\geq\lambda_1(\Om^*)\ \ \ \mbox{with equality if and only if}\ \ \Omega\ \ \mbox{is a ball,}
\end{eqnarray}
which was proved by Faber and Krahn independently in \cite{Fab} and \cite{Kra} respectively. The second result we recall is the following Ashbaugh-Benguria inequality
\begin{eqnarray}\label{eq1.3}
\frac{\la_2(\Om)}{\la_1(\Om)}\leq\frac{\la_2(\Omega^*)}{\la_1(\Omega^*)}\ \ \ \mbox{with equality if and only if}\ \ \Omega\ \ \mbox{is a ball,}
\end{eqnarray}
which is a conjecture of Payne, P$\acute{o}$lya  and Weinberger in \cite{PPW1,PPW2} for dimension $n=2$, and of Thompson in \cite{Thompson} for dimension $n>2$. Eventually, this conjecture was proved by M.S.Ashbaugh and R.D.Benguria in \cite{Ben} for dimension $n=2$, and in \cite{Ashb} for dimension $n>2$. The last result we recall here is the following Payne-Rayner inequality for the first eigenfunction $\varphi_1(x)$ of problem (\ref{eq12.4}) in dimension $n=2$.
\begin{eqnarray}\label{eq1.4}
\int_{\Omega}\varphi_1^2(x)dx\leq\frac{\lambda_1(\Omega)}{4\pi}(\int_{\Omega}\varphi_1(x)dx)^2
\end{eqnarray}
with equality if and only if $\Omega$ is a disk. The above inequality was proved by L.E.Payne and M.E.Rayner in \cite{PayRay}, and successively generalized to any dimension by M.Th$\acute{e}$re$\grave{s}$e and K.Jobin in \cite{Jobin} and by G.Chiti in \cite{ref GC} with method different from that of \cite{PayRay}. It is worth pointing out that G.Chiti has in fact proved a reverse H$\ddot{o}$lder inequality in \cite{ref GC} which is more general than the Payne-Rayner inequality.

When $0<\be<+\infty$, problem (\ref{eq1.1}) is called eigenvalue problem for Robin Laplacian. There are also some results for the eigenvalue problem of Robin Laplacian though it is few. At first, for any $\beta>0$, we have the following Faber-Krahn type inequality
\begin{eqnarray}\label{eq12.6}
\lambda_1(\Om,\be)\geq\lambda_1(\Om^*,\be)\ \ \ \mbox{with equality if and only if}\ \ \Omega\ \ \mbox{is a ball,}
\end{eqnarray}
which was proved by Bossel in \cite{Bos} for dimension $n=2$, and by Danners in \cite{Dane1,Dane2} for dimension $n>2$. It is worthy of mention that inequality (\ref{eq12.6}) was recently generalized by Q.Y.Dai and Y.X.Fu in \cite{QD} to the Robin problem involving p-Laplacian. In the second, Payne and Schaefer proved the following estimate for the ratio of the first two eigenvalues in \cite{payne}.
\begin{eqnarray}\label{eq12.7}
\frac{\la_2(\Om,\be)}{\la_1(\Om,\be)}\leq1+\frac{4}{n}\ \ \ \mbox{for}\ \ \be>P_0\la_1(\Om)
\end{eqnarray}
with $P_0=\max\limits_{x\in\partial\Om}x\cdot\nu$. The inequality (\ref{eq12.7}) is an extension of Payne, P$\acute{o}$lya  and Weinbergers result in \cite{PPW1,PPW2}, and of Thompson's result in \cite{Thompson}. Obviously, inequality (\ref{eq12.7}) can not be valid for all $\beta>0$ since
\begin{eqnarray*}
\frac{\la_2(\Om,\be)}{\la_1(\Om,\be)}\rightarrow +\infty\ \ \ \mbox{as}\ \ \ \beta\rightarrow 0^+
\end{eqnarray*}
due to (\ref{eq12.2}).

Motivated by the inequality (\ref{eq12.7}), A. Henrot proposed a question that for what $\be$ the ratio $\frac{\la_2(\Om,\be)}{\la_1(\Om,\be)}$ achieves its maximum for the ball in a recent paper \cite{Henrot1}. From (\ref{eq12.3}), (\ref{eq1.3}) and (\ref{eq12.6}), one can see that the answer to the Henrot's problem should be positive for the parameter $\be$ small, or large enough. This leads us to make a conjecture as the following




{\em{\bf Conjecture A.} For any $\be>0$, there holds
\begin{eqnarray}\label{eq12.8}
\frac{\la_2(\Om,\be)}{\la_1(\Om,\be)}\leq\frac{\la_2(\Om^*,\be)}
{\la_1(\Om^*,\be)}
\end{eqnarray}
and the equality occurs if and only if $\Om$ is a ball.}

At last, we point out here that the Payne-Rayner inequality was also partially generalized by F.Takahashi and A.Uegaki in \cite{Takahashi} from Dirichlet Laplacian to Robin Laplacian (see also \cite{QWang} for more information).

The aims of this paper are two folds. One is to shed some lights on the proof of conjecture A; the other is to extend the Chiti's reverse H$\ddot{o}$lder inequality to the first eigenfunction of problem (\ref{eq1.1}) with parameter $\beta\in(0,+\infty)$. To this end, a crucial step is to prove a Chiti type comparison result for problem (\ref{eq1.1}) with $\beta\in(0,+\infty)$.

Let $|\Omega|$ denote the volume of domain $\Omega$, and $\omega_n$ be the volume of the unit ball in $R^n$. Set
$$R^*=\left(|\Om|/\om_n\right)^{\frac{1}{n}}\ \ \ \mbox{and}\ \ \ \rho=(\sqrt{\lambda_1(\Om^*)}/\sqrt{\lambda_1(\Om)})R^*.$$
If we denote by $B_\rho(0)$ the ball in $R^n$ with radius $\rho$ and center at origin, and by $Y_1(x)$ the first eigenfunction of the eigenvalue problem
\begin{eqnarray}
\left\{\begin{array}{ll}
-\Delta Y=\lambda Y & x\in B_\rho(0),\\
Y=0 & x\in\partial B_\rho(0),
\end{array}\label{eq1.5}
\right.
\end{eqnarray}
then the Chiti's comparison result for Dirichlet Laplacian, that is, for problem (\ref{eq12.4}) can be stated as the following



{\em{\bf Theorem B(\cite{ref GC}).}  Let $Y^*_1(s)$ and $\varphi^*_1(s)$ are the decreasing rearrangement of $Y_1$ and $\varphi_1$, whose definition is given in section 2, respectively. If, for $p>0$, we normalize $\varphi_1(x)$ and $Y_1(x)$ so that $\int_{\Om}\varphi^p_1dx=\int_{B_\rho(0)}Y^p_1dx$, then there exists a unique point $s_0\in (0,\ |B_\rho(0)|)$ such that
\begin{eqnarray}
\left\{\begin{array}{ll}
Y^*_1(s)>\varphi^*_1(s)  & \mbox{for\ }s\in (0,\ s_0),\\
Y^*_1(s)\leq\varphi^*_1(s) & \mbox{for\ }s\in [s_0,\ |B_\rho(0)|),
\end{array}\label{eq1.6}
\right.
\end{eqnarray}}

Chiti's comparison result was proved by making use of the Schwarz symmetrization method. This method requires an application of the classical isoperimetric inequality to the level set $\{x\in\Om:\ \varphi_1(x)>t\}$ of $\varphi_1(x)$. It is well known that the classical isoperimetric inequality can only be used in the case where the boundary of the domain under consideration is a closed surface. Hence, Chiti can prove his comparison result fortunately on the full interval $(0,\ |B_\rho(0)|)$ due to the fact that the boundary of the level set $\{x\in\Om:\ \varphi_1(x)>t\}$ of the first Dirichlet eigenfunction $\varphi_1(x)$ is indeed a closed surface for any $t>0$. However, the level surface $\{x\in\Omega:\ \psi_1(x)=t\}$ of the first Robin eigenfunction $\psi_1(x)$ is always not a closed surface for $t>0$ small enough provided that $\Omega$ is not a ball. Hence, we can not expect to establish a Chiti type comparison result, which is good enough to solve conjecture A completely, for Robin problem (\ref{eq1.1}). This may be the essential difficulty in the study of conjecture A. The main observation of this paper is that the level surface $\{x\in\Omega:\ \psi_1(x)=t\}$ is a closed surface for $t$ large in some extent, which can be used to establish a Chiti type comparison result for Robin problem (\ref{eq1.1}) on a small interval. Once the Chiti type comparison result is established, we can follow the arguments used in \cite{Ben}, \cite{Ashb} and \cite{ref GC} to get main results of this paper.


To state our results precisely, we fix some notations first. We always assume that $\beta\in(0,+\infty)$, and $\psi_1(x)$ is the first eigenfunction of problem (\ref{eq1.1}) in the following paragraphs. Let
$$M=\max\limits_{x\in\partial\Omega}\psi_1(x)\ \ \ \mbox{and}\ \ \ \Om_M=\{x\in\Om:\ \psi_1(x)>M\}.$$
It is easy to see that the boundary $\{x\in\overline{\Om}:\ \psi_1(x)=M\}$ of $\Om_M$ is a closed surface. Furthermore, for any $t>M$, the level surface $\{x\in\Om:\ \psi_1(x)=t\}$ of $\psi_1(x)$ is also a closed surface. Hence, the classical isoperimetric inequality can be applied to any level set $\Om_t=\{x\in\Om:\ \psi_1(x)>t\}$ of $\psi_1(x)$ when $t\geq M$.

Let
$$R_\lambda=(\sqrt{\lambda_1(\Om^*,\be)}/\sqrt{\lambda_1(\Om,\be)})R^*.$$
Then, from the dilation of problem (\ref{eq1.1}) and the inequality (\ref{eq12.6}), we have
\begin{eqnarray}\label{eq12.9}
\lambda_1(\Omega,\beta)=\lambda_1(B_{R_\lambda},\frac{R^*}{R_\lambda}\beta)\ \ \ \mbox{and}\ \ \  R_\lambda\leq R^*.
\end{eqnarray}

Let
\begin{eqnarray}\label{eq12.10}
R_M=(|\Om_M|/\omega_n)^{\frac{1}{n}}\ \ \ \mbox{and}\ \ \ R=\mbox{min\ }\{R_\lambda,R_M\}.
\end{eqnarray}
It is easy to see that $R_M$ depends only on $\beta$, $\Omega$ and $n$ since the first eigenfunction $\psi_1(x)$ is unique up to multiplication of a positive constant (see \cite{QD}), and the set $\Om_M=\{x\in\Om:\ \psi_1(x)>M\}$ is independent of the choice of $\psi_1(x)$. Denote by $\lambda_1(B_R,\frac{R^*}{R}\beta)$ the first eigenvalue, and $z_1(x)$ the first eigenfunction of the
following eigenvalue problem
\begin{eqnarray}\label{eq1.7}
\left\{
\begin{array}{ll}
-\Delta z=\lambda z &x\in B_R(0),\\
\frac{\partial z}{\partial\nu}+\frac{R^*}{R}\beta z=0 &x\in\partial B_R(0).
\end{array}
\right.
\end{eqnarray}
By (\ref{eq12.9}), (\ref{eq12.10}) and a result of T.Giorgi and R.G.Smits in \cite{Gio}, we always have

\begin{eqnarray}\label{eq12.11}
\lambda_1(B_R,\frac{R^*}{R}\beta)\geq\lambda_1(\Omega,\beta)
\end{eqnarray}

Keeping all above notations in mind, the first result of our paper can be stated as

\vskip 0.1cm

{\em{\bf Theorem 1.1.}\ For any $\beta>0$, we have the following estimate
\begin{equation}\label{eq1.8.0}
  \frac{\la_2(\Om,\be)}{\la_1(\Om,\be)}\leq\frac{R_\lambda^2}{R^2}\frac{\la_2(\Om^*,\be)}{\la_1(\Om^*,\be)}-\frac{R_\lambda^2}{R^2}+1.
\end{equation}
}
\vskip 0.1cm
{\em{\bf Remark 1.2.}\ If $\Omega$ is a ball, then equality occurs in (\ref{eq1.8.0}). In fact, by a result of Q.Y.Dai and Y.X.Fu in \cite{QD}, we know that the first eigenfunction for Robin Laplacian on a ball is radially symmetry and decreasing. Combining this observation with the Faber-Krahn type inequality (\ref{eq12.6}), we can see that $R_M=R_\lambda=R^*$. Hence, $R_\lambda^2/R^2=1$,  and we get the equality in (\ref{eq1.8.0}).}
\vskip 0.1cm

{\em{\bf Corollary 1.3.}\ If $R_M\geq R_\lambda$, then $R=R_\lambda$, and the inequality (\ref{eq1.8.0}) becomes
  $$ \frac{\la_2(\Om,\be)}{\la_1(\Om,\be)}\leq\frac{\la_2(\Om^*,\be)}{\la_1(\Om^*,\be)}.$$}

\vskip 0.1cm
{\em{\bf Remark 1.4.}\ If $\beta=+\infty$,  we have $M=0$ and $\Omega_M=\Omega$. Thus, $R_M=R^*\geq R_\lambda$, and the Ashbaugh-Benguria inequality can be recovered from the conclusion of corollary 1.3.}

\vskip 0.1cm
{\em{\bf Remark 1.5.}\ Though, the exact value of $M$ and $R_M$ is not known for general domain $\Omega$, we can get the following rough estimate of $R_M$ for convex domains in section 5.
\begin{equation}\label{eq1.8.1}
R_M\geq\bigg[\frac{2n}{\lambda_1(\Omega)}\bigg(1-\sqrt{\frac{\frac{2}{n}\la_1(\Om)}
{\be^2+\frac{2}{n}\la_1(\Om)}}\bigg)\bigg]^{\frac{1}{2}}.
\end{equation}
}
\vskip 0.2cm

The second result of our paper is the following Chiti type reverse H$\ddot{o}$lder inequality

\vskip 0.1cm
 {\em{\bf Theorem 1.6.}
 For any $q\geq p>0$, there holds
  \begin{eqnarray*}
\left(\int_{\Om}\psi^q_1dx\right)^{\frac{1}{q}}\leq  K\left(p,q,\beta,\Omega,n\right)\left(\int_{\Om}\psi^p_1dx\right)^{\frac{1}{p}},
\end{eqnarray*}
 where $K\left(p,q,\beta,\Omega,n\right)$ is a positive constant will be given in section 4.}
\vskip 0.1cm
The rest part of this paper is organized as follows: Section 2 is a collection of some basic facts about the rearrangement of nonnegative measurable functions. Section 3 includes a proof of Chiti type comparison result. The proofs of Theorem 1.1 and 1.6 are presented in Section 4. A detailed explanation of Remark 1.5 is given in Section 5. An appendix is arranged to give some Lemmas needed in the proofs of Theorem 1.1 and 1.6.
\vskip 0.1in
\section*{2.\ Preliminary}

\setcounter{section}{2}

\setcounter{equation}{0}

\renewcommand{\theequation}{\thesection.\arabic{equation}}

\vskip 0.1cm

  In this section, we recall some basic facts about the rearrangement
of nonnegative measurable functions.

Let $f:\ \Om\mapsto R$ be a nonnegative measurable function. For any
$t\geq 0$. The level set $\Om_t$ of $f$ at the level $t$ is defined
by
$$\Om_t=\{x\in \Om\ |\ f(x)>t\},\ \ \ \ t\geq 0.$$
The distribution function of $f$ is given by
$$\mu_f(t)=|\Om_t|=\mbox{meas}\{x\in \Om\ |\ f(x)>t\},\ \ \ \ t\geq 0.$$
Obviously, $\mu_f(t)$ is a monotonically decreasing function of $t$, $\mu_f(t)=0$ for $t\geq \mbox{ess}\sup f(x)$, and
$\mu_f(t)=|\Om|$ for $t=0$.

\vskip 0.1cm

{\bf Definition 2.1.}  Let $\Om$ be a bounded domain in $R^n$, $f:\
\Om\mapsto R$ be a nonnegative measurable function. Then the
decreasing rearrangement $f^*$ of $f$ is a function defined on $[0,\
\infty)$ by
\begin{eqnarray*}
f^*(s)=\left
 \{\begin{array}{ll}
 \mbox{ess}\sup\limits_{x\in\Om}f(x) &\mbox{for\ }s=0,\\
 \\
 \inf\{t>0|\mu_f(t)<s\} &\mbox{for\ }s>0.
\end{array}
\right.
\end{eqnarray*}
Obviously, $f^*(s)=0$, for $s\geq |\Om|$. The increasing rearrangement $f_*$ of $f$ is defined by $f_*(s)=f^*(|\Om|-s)$ for $s\in(0,\ +\infty)$.

\vskip 0.1cm

{\bf Definition 2.2.} Let $\Om$ be a bounded domain in $R^n$, $f:\
\Om\mapsto R$ be a nonnegative measurable function. Then the
decreasing Schwarz symmetrization $f^{\star}$ of $f$ is a function
defined by
$$f^{\star}(x)=f^{*}(\om_n|x|^n),\ \ \ \mbox{for\ }x\in\Om^*.$$

There are many fine properties of rearrangement. Here we only collect some important properties needed in this paper.

\vskip 0.1cm

{\bf Proposition 2.3.} Let $f:\ \Om\mapsto R$ be a
nonnegative measurable function. Then, $f,\ f^*$ and $f^{\star}$ are
all equimeasurable and
$$\int_{\Om}fdx=\int^{|\Om|}_0f^*(s)ds=\int_{\Om^*}f^{\star}(x)dx.$$
Moreover,  for any Borel measurable function $F:\ R\mapsto R$, there
holds

$$
\int_{\Om}F(f(x))dx=\int^{|\Om|}_0F(f^*(s))ds
=\int_{B_{R^*}(0)}F(f^{\star}(x))dx.$$

\vskip 0.1cm

{\bf Proposition 2.4.} If $f:\ [0,\ l]\mapsto R$ is
nonnegative and non-increasing, then $f=f^* \ a.e.$

\vskip 0.1cm
{\bf Proposition 2.5.} If $\psi:\ R\mapsto R$ is a non-decreasing function, then
$$\psi(f^*)=(\psi(f))^*,\ \ \psi(f^{\star})=(\psi(f))^{\star}$$
for any nonnegative measurable function  $f:\ \Om\mapsto R$.

\vskip 0.1cm
{\bf Proposition 2.6.} Let $f\in L^p(\Om),\ g\in L^q(\Om)$ with $\frac{1}{p}+\frac{1}{q}=1$. Then
$$
\int^{|\Om|}_0f^*(s)g_*(s)ds\leq\int_{\Om}f(x)g(x)dx\leq\int^{|\Om|}_0f^*(s)g^*(s)ds,
$$
$$
\int_{\Om^*}f^{\star}(x)g_{\star}(x)dx\leq\int_{\Om}f(x)g(x)dx\leq\int_{\Om^*}f^{\star}(x)g^{\star}(x)dx.
$$
Consequently
$$
\int_{E}f(x)dx\leq\int^{|E|}_0f^*(s)ds=\int_{E^*}f^{\star}(x)dx.
$$
for any measurable set $E\subset\Om$.

 \vskip 0.1cm

{\bf Proposition 2.7.} If $f(x)=f(|x|)$ is nonnegative, and is decreasing (or increasing) as a function of
$r=|x|$ for $x\in\Om$, then
$$
f^{\star}(r)\leq f(r) (\mbox{or}\ f_{\star}(r)\geq f(r))\ \ \ \mbox{for\ }r\in(0,R^*).
$$
\vskip 0.1cm

{\bf Proposition 2.8.} Let $T,\ \alpha,\ \beta$ be
real numbers such that $0<\alpha\leq \beta$ and $T>0$. If $f,\ g$
are real functions in $L^\beta([0,\ T])$, then we have
$$\int^{T}_0f^{*^{\beta}}(t)dt\leq\int^{T}_0g^{*^{\beta}}(t)dt.$$
provided that
$$\int^{s}_0f^{*^{\alpha}}(t)dt\leq\int^{s}_0g^{*{^{\alpha}}}(t)dt \ \ \ \mbox{for any}\ \ \ s\in[0,\ T].$$

\vskip 0.1in
For detailed information of all the above propositions, we refer to \cite{ref GJG}, \cite{ref BK} and \cite{ref SK}.

\section*{3.\ Chiti Type Comparison Result}

\setcounter{section}{3}

\setcounter{equation}{0}
This section devotes to prove a Chiti type comparison result
for problem (\ref{eq1.1}). Keeping notations given in section 1 in use, our Chiti type comparison result can be stated as
\vskip 0.1cm
{\em{\bf Theorem 3.1.}\ For any $p>0$, if we normalize $\psi_1(x)$ so that $\int_{\Om}\psi^p_1dx=\int_{B_R(0)}z^p_1dx$, then the following statements hold.
 \vskip 0.05cm
(i)\ In the case $z^*_1(|B_R(0)|)\geq\psi^*_1(|B_R(0)|)$, $z^*_1(s)\geq\psi^*_1(s)$ for any $s\in (0,\ |B_R(0)|)$.
 \vskip 0.05cm
(ii)\ In the case $z^*_1(|B_R(0)|)<\psi^*_1(|B_R(0)|)$, there exists a unique $s_0\in (0,\ |B_R(0)|)$ such that

\begin{eqnarray*}
\left\{\begin{array}{ll}
z^*_1(s)\geq\psi^*_1(s)  & \mbox{for\ }s\in [0,\ s_0],\\
z^*_1(s)<\psi^*_1(s) & \mbox{for\ }s\in(s_0, |B_R(0)|].
\end{array}
\right.
\end{eqnarray*}}

\vskip 0.1cm


The proof of Theorem 3.1 depends strongly on the following lemma. Hence, we stop to give a proof of it before proceeding on.

{\em {\bf Lemma 3.2.} Assume that $0<\beta< +\infty$. Then the following inequality holds for any $s\in(0,|\Omega_M|)$.
\begin{eqnarray}
-\frac{d\psi^*_1(s)}{ds} \leq
n^{-2}w_n^{-\frac{2}{n}}s^{\frac{2}{n}-2}
\int^s_0\lambda_1(\Omega,\beta)\psi^*_1(\tau)d\tau.
\label{eq1.11}
\end{eqnarray}
}
{\bf Proof:}\ Since $\psi_1(x)$ satisfies
\begin{eqnarray*}
-\Delta \psi_1(x)=\lambda_1(\Omega,\beta)\psi_1(x) & x\in\Omega,\\
\end{eqnarray*}
we have
\begin{eqnarray*}
-\int_{\Omega_t}\Delta
\psi_1(x)dx=\lambda_1(\Omega,\beta)\int_{\Omega_t}\psi_1(x)dx.
\end{eqnarray*}
Noticing that $\Omega_t\subset\subset\Omega$ for any $t>M$, we have
 \begin{eqnarray*}
-\int_{\Omega_t}\Delta \psi_1(x)dx=-\int_{\partial\Omega_t}
\frac{\partial\psi_1(x)}{\partial\nu}dS=\int_{\partial\Omega_t}|\nabla\psi_1(x)|dS,
\end{eqnarray*}
and
\begin{equation*}
\int_{\partial\Omega_t}|\nabla\psi_1(x)|dS
\int_{\partial\Omega_t}\frac{1}{|\nabla\psi_1(x)|}dS
\geq|\partial\Omega_t|^2,
\end{equation*}
Hence,
\begin{eqnarray*}
-\int_{\Omega_t}\Delta \psi_1(x)dx\geq\frac{|\partial\Omega_t|^2}
{\int_{\partial\Omega_t}\frac{1}{|\nabla\psi_1(x)|}dS}.
\end{eqnarray*}
By the co-area formula, we have
\begin{eqnarray*}
\mu_{\psi_1}(t)=|\Omega_t|=\int_{\Omega_t}dx=\int^{+\infty}_{t}
{\int_{\partial\Omega_\tau}\frac{dS}{|\nabla\psi_1(x)|}d\tau}.
\end{eqnarray*}
Consequently,
\begin{eqnarray*}
{{\mu'_{\psi_1}}(t)}=\frac{d\mu_{\psi_1}(t)}{dt}=-\int_{\partial\Omega_t}\frac{dS}{|\nabla\psi_1(x)|},
\end{eqnarray*}
and
\begin{eqnarray*}
-\int_{\Omega_t}\Delta \psi_1(x)dx\geq-\frac{|\partial\Omega_t|^2}{{\mu'_{\psi_1}}(t)}.
\end{eqnarray*}
Since $\partial\Omega_t$ is a closed surface when $t\geq M$, we can apply the classical isoperimetric inequality to get
\begin{eqnarray*}
|\partial\Omega_t|\geq nw_n^\frac{1}{n}|\Omega_t|^{1-\frac{1}{n}}=nw_n^\frac{1}{n}\mu_{\psi_1}^{1-\frac{1}{n}}\ \ \ \mbox{for}\ \ \ t\geq M.
\end{eqnarray*}
This implies that
\begin{eqnarray*}
-\int_{\Omega_t}\Delta \psi_1(x)dx\geq\frac{n^2w_n^\frac{2}{n}\mu_{\psi_1}^{2-\frac{2}{n}}}{-{\mu'_{\psi_1}}(t)}\ \ \ \mbox{for}\ \ \ t\geq M.
\end{eqnarray*}
Noting that
\begin{eqnarray*}
\lambda_1(\Omega,\beta)\int_{\Omega_t}\psi_1(x) dx =\lambda_1(\Omega,\beta)
\int^{\mu_{\psi_1}}_{0}\psi^*_1(\tau)d\tau,
 \end{eqnarray*}
we get
\begin{eqnarray*}
-\frac{1}{{\mu'_{\psi_1}}(t)}\leq\lambda_1(\Omega,\beta)n^{-2}w_n^{-\frac{2}{n}}\mu_{\psi_1}^{\frac{2}{n}-2}
\int^{\mu_{\psi_1}}_0\psi^*_1(\tau)d\tau.
\end{eqnarray*}
Since $\psi^*_1(s)$ is essentially a reverse function of $\mu_{\psi_1}(t)$, we have
\begin{eqnarray*}
-\frac{d\psi^*_1(s)}{ds} \leq n^{-2}w_n^{-\frac{2}{n}}s^{\frac{2}{n}-2}
\int^s_0\lambda_1(\Omega,\beta)\psi^*_1(\tau)d\tau\ \ \ \mbox{for}\ \ \ s\in(0,|\Omega_M|).
\end{eqnarray*}
This is just the desired conclusion of Lemma 3.2.

\vskip 0.1cm
{\em {\bf The proof of Theorem 3.1:}}\ From Lemma 3.2, (\ref{eq12.10}) and (\ref{eq12.11}), we know that $\psi^*_1(s)$ satisfies
\begin{eqnarray}\label{eq1.26.1}
-\frac{d\psi^*_1(s)}{ds} \leq n^{-2}w_n^{-\frac{2}{n}}s^{\frac{2}{n}-2}
\int^s_0\lambda_1(B_R,\frac{R^*}{R}\beta)\psi^*_1(\tau)d\tau\ \ \ \mbox{for}\ \ \ s\in(0,|B_R(0)|).
\end{eqnarray}
By Proposition 2.4, Proposition 2.5 and (\ref{eq1.7}), we deduce that $z^*_1(s)$ satisfies
\begin{eqnarray}\label{eq1.26.2}
-\frac{dz^*_1(s)}{ds} = n^{-2}w_n^{-\frac{2}{n}}s^{\frac{2}{n}-2}
\int^s_0\lambda_1(B_R,\frac{R^*}{R}\beta)z^*_1(\tau)d\tau\ \ \ \mbox{for}\ \ \ s\in(0,|B_R(0)|).
\end{eqnarray}
At this stage, we divide the proof of Theorem 3.1 into two cases.

(i) In the case $z^*_1(|B_R(0)|)\geq\psi^*_1(|B_R(0)|)$, we want to prove $z^*_1(s)\geq\psi^*_1(s)$ for any $s\in (0,\ |B_R(0)|)$.
If this conclusion is not true, then there should exist an interval
$(s_1,s_2) \subset(0,|B_R(0)|)$ such that
$z^*_1(r)<\psi^*_1(r)$ for $r\in(s_1,s_2)$, and $z^*_1(s_i)=\psi^*_1(s_i)$  for $i=1,2$. It follows from the assumption $\int_{\Om}\psi^p_1dx=\int_{B_R(0)}z^p_1dx$ that either $s_1\neq 0$, or $s_2\neq |B_R(0)|$. No loss of generality, we assume that $s_2\neq |B_R(0)|$. Choosing
$$s_2=\inf\{s:\ z^*_1(\tau)\geq\psi^*_1(\tau),\tau\in(s,|B_R(0)|)\},$$
it is easy to see that $s_2\neq 0$ and $z^*_1(s_2)=\psi^*_1(s_2)$. Fixing $s_2$, we choose
$$s_1=\inf\{s:\ z^*_1(\tau)<\psi^*_1(\tau), \tau\in(s,s_2)\}.$$
Then, there are two possibilities for $s_1$. One is $s_1=0$, and the other is $s_1\neq 0$.

If $s_1=0$, or $s_1\neq 0$ and $\int^{s_1}_0\psi^*_1(\tau)d\tau>\int^{s_1}_0 z^*_1(\tau)d\tau$, we let
\begin{equation*}
w(s)=\left\{\begin{array}{ll}
\psi^*_1(s),\ \ \ s\in[0,s_2)\\
z^*_1(s), \ \ \ s\in[s_2,|B_R(0)|).
\end{array}
\right.
\end{equation*}

If $s_1\neq0$, and $\int^{s_1}_0\psi^*_1(\tau)d\tau\leq\int^{s_1}_0 z^*_1(\tau)d\tau$, we let
\begin{equation*}
w(s)=\left\{\begin{array}{ll}
z^*_1(s),\ \ \ s\in[0,s_1]\\
\psi^*_1(s),\ \ \ s\in[s_1 ,s_2],\\
z^*_1(s),\ \ \ s\in[s_2,|B_R(0)|].
\end{array}
\right.
\end{equation*}

It is easy to check that $w(s)$ satisfies

\begin{eqnarray}\label{eq3.1}
-\frac{dw(s)}{ds} \leq n^{-2}w_n^{-\frac{2}{n}}s^{\frac{2}{n}-2}\int^s_0\lambda_1(B_R,\frac{R^*}{R}\beta)w(\tau)d\tau\ \ \ \mbox{for any}\ \ \ s\in (0,|B_R(0)|).
\end{eqnarray}

Define a test function $W(x)$ of $\lambda_1(B_R,\frac{R^*}{R}\beta)$ by  $W(x)=w(w_n|x|^n)$ for any $x\in B_R(0)$. By the definition of $w(s)$, we see that $W(x)\not\equiv z_1(x)$. Hence, we have

\begin{equation}
\begin{array}{ll}
\lambda_1(B_R,\frac{R^*}{R}\beta)\int_{B_R}W^2 dx
& <\int_{B_R}|\nabla W|^2dx+\frac{R^*\beta}{R}\int_{\partial B_R}W^2 dS\\
& =\int_{B_R}|\nabla W|^2dx+\frac{R^*\beta}{R}z^2_1(R)|\partial B_R|.
\end{array}\label{eq3.0}
\end{equation}
Since
\begin{eqnarray*}
\int_{B_R}|\nabla W|^2dx=\int^{|B_R|}_0({w'(s)})^2 n^{2}w_n^{\frac{2}{n}}s^{2-\frac{2}{n}}ds,
\end{eqnarray*}
it follows from (\ref{eq3.1}) that
\begin{eqnarray}\label{eq3.2}
\int_{B_R}|\nabla W|^2dx\leq -\int^{|B_R|}_0{w'(s)}\left(\lambda_1(B_R,\frac{R^*}{R}\beta)\int^s_0w(\tau)d\tau \right)ds.
\end{eqnarray}
Let
$$I=-\int^{|B_R|}_0{w'(s)}\left(\lambda_1(B_R,\frac{R^*}{R}\beta)\int^s_0w(\tau)d\tau \right)ds.$$
By integration by parts, we can get
\begin{eqnarray}\label{eq3.3}
I= \lambda_1(B_R,\frac{R^*}{R}\beta)\bigg(\int_{B_R}W^2(x)dx-
                        w(|B_R|)\int^{|B_R|}_0w(\tau)d\tau\bigg).
\end{eqnarray}
Setting
$$II=\lambda_1(B_R,\frac{R^*}{R}\beta)w(|B_R|)\int^{|B_R|}_0w(\tau)d\tau,$$
it follows from (\ref{eq3.2}) and (\ref{eq3.3}) that
\begin{eqnarray}\label{eq3.4}
\int_{B_R}|\nabla W|^2dx\leq\lambda_1(B_R,\frac{R^*}{R}\beta)\int_{B_R}W^2(x)dx-II.
\end{eqnarray}

If $s_1=0$, or $s_1\neq 0$ and $\int^{s_1}_0\psi^*_1(\tau)d\tau>\int^{s_1}_0 z^*_1(\tau)d\tau$, by virtue of (\ref{eq1.26.1}), (\ref{eq1.26.2}) and the definition of $w(s)$, we can estimate $II$ as
\begin{equation}\label{eq3.5}
\begin{array}{ll}
II &=\lambda_1(B_R,\frac{R^*}{R}\beta)z_1(R)\left[\int^{s_2}_0\psi^*_1(\tau)d\tau+\int^{|B_R|}_0z^*_1(\tau)d\tau- \int^{s_2}_0z^*_1(\tau)d\tau\right]\\
&\geq z_1(R)n^{2}w_n^{\frac{2}{n}}\left[s_2^{2-\frac{2}{n}}\frac{d(z^*_1-\psi^*_1)(s_2)}{ds}
                                    - |B_R|^{2-\frac{2}{n}}\frac{dz^*_1(|B_R|)}{ds}\right].
\end{array}
\end{equation}
Since $s=\omega_nr^n=\omega_n|x|^n$ for $x\in B_R(0)$, we have
$$\frac{dz^*_1(|B_R|)}{ds}=\frac{dz_1(r)}{dr}
\frac{dr}{ds}\bigg|_{r=R}=n^{-1}w_n^{-1}R^{1-n}\frac{dz_1(R)}{dr}.$$
By the boundary condition of $z_1(x)$, we have
$$\frac{dz_1(R)}{dr}=-\frac{R^*\beta}{R}z_1(R).$$
Therefore
\begin{eqnarray}\label{eq3.6}
\frac{dz^*_1(|B_R|)}{ds}=-\frac{R^*\beta z_1(R)}{n|B_R(0)|}
\end{eqnarray}
Substituting (\ref{eq3.6}) into (\ref{eq3.5}), we get
\begin{eqnarray}\label{eq3.7}
 II\geq n^{2}w_n^{\frac{2}{n}}s_2^{2-\frac{2}{n}}\frac{d(z^*_1-\psi^*_1)(s_2)}{ds}z_1(R)
                                    +\frac{R^*\beta}{R}z^2_1(R)|\partial B_R|.
\end{eqnarray}
Combining (\ref{eq3.0}), (\ref{eq3.4}) with (\ref{eq3.7}), we reach
\begin{eqnarray}\label{eq3.8}
\frac{d(\psi^*_1-z^*_1)}{ds}(s_2)>0.
\end{eqnarray}
Since $\psi^*_1(s_2)-z^*_1(s_2)=0$ and $\psi^*_1(s)-z^*_1(s)>0$ for
$s\in(s_1,\ s_2)$, we have
 $$\frac{d(\psi^*_1-z^*_1)}{ds}(s_2)\leq0.$$
This contradicts (\ref{eq3.8}).

If $s_1\neq 0$, and $\int^{s_1}_0\psi^*_1(\tau)d\tau\leq\int^{s_1}_0 z^*_1(\tau)d\tau$, by virtue of (\ref{eq1.26.1}), (\ref{eq1.26.2}) and the definition of $w(s)$, we can estimate $II$ as
\begin{eqnarray}\label{eq3.9}
II=\lambda_1(B_R,\frac{R^*}{R}\beta)z_1(R)(\int^{s_1}_0z^*_1(\tau)+\int^{s_2}_{s_1}\psi^*_1(\tau)
+\int^{|B_R|}_{s_2}z^*_1(\tau)).
\end{eqnarray}
Since $\psi^*_1(\tau)>z^*_1(\tau)$ for $\tau\in(s_1,\ s_2)$, we have
\begin{eqnarray}\label{eq3.10}
\int^{s_2}_{s_1}\psi^*_1(\tau)\geq\int^{s_2}_{s_1}z^*_1(\tau).
\end{eqnarray}
Substituting (\ref{eq3.10}) into (\ref{eq3.9}), we get
\begin{equation}
\begin{array}{ll}
II & \geq\lambda_1(B_R,\frac{R^*}{R}\beta)z_1(R)(\int^{s_1}_0z^*_1(\tau)+\int^{s_2}_{s_1}z^*_1(\tau)+\int^{|B_R|}_{s_2}z^*_1(\tau))\\
   &=-n^{2}\omega_n^2R^{2n-2}z_1(R)\frac{dz^*_1(|B_R|)}{ds}\\
   &=-|\partial B_R|z_1(R)\frac{dz_1(R)}{dr}.
\end{array}\label{eq3.11}
\end{equation}
Putting the boundary condition into (\ref{eq3.11}), we arrive
\begin{eqnarray}\label{eq3.12}
II\geq\frac{R^*}{R}\beta z_1^2(R)|\partial B_R|.
\end{eqnarray}
Inserting (\ref{eq3.12}) and (\ref{eq3.4}) into (\ref{eq3.0}), we deduce that
$$\lambda_1(B_R,\frac{R^*}{R}\beta)\int_{B_R}W^2 dx<\lambda_1(B_R,\frac{R^*}{R}\beta)\int_{B_R}W^2 dx.$$
A contradiction. This completes the proof of Theorem 3.1 $(i)$.

(ii) In the case $z^*_1(|B_R|)<\psi^*_1(|B_R|)$, we want to prove there exists a unique point $s_0\in(0,|B_R(0)|)$ such that
\begin{eqnarray*}
\left\{\begin{array}{ll}
z^*_1(s)\geq\psi^*_1(s)  & \mbox{for}\ \ s\in [0,\ s_0],\\
z^*_1(s)<\psi^*_1(s) & \mbox{for}\ \ s\in(s_0, |B_R(0)|].
\end{array}
\right.
\end{eqnarray*}
At first, from the assumptions that
$$z^*_1(|B_R|)<\psi^*_1(|B_R|)\ \ \ \mbox{and}\ \ \ \int^{|\Om^*|}_0(\psi^*_1)^p(s)ds=\int^{|B_R(0)|}_0(z^*_1)^p(s)ds,$$
we can easily see that $\psi^*_1(s)$ and $z^*_1(s)$ must intersect at some point $s_0\in(0,\ |B_R|)$.
Choosing
$$s_0=\inf\{s:\ z^*_1(\tau)<\psi^*_1(\tau),\tau\in(s,|B_R(0)|)\},$$
we are going to prove that $s_0$ is the unique point we want. If this is not true, we can find a point $s_1\in(0,s_0)$ such that
$$z^*_1(s)\geq\psi^*_1(s)\ \ \ \mbox{for any}\ \ \ s\in(s_1,s_0),\ \ \ \mbox{and}\ \ \ z^*_1(s)\not\equiv\psi^*_1(s)\ \ \ \mbox{on}\ \ \ (s_1,s_0)$$
due to the assumption $\int^{|\Om^*|}_0(\psi^*_1)^p(s)ds=\int^{|B_R(0)|}_0(z^*_1)^p(s)ds$.

Let
\begin{equation*}
w(s)=\left\{\begin{array}{ll}
\psi^*_1(s), & s\in[0,s_1],\ \mbox{if}\ \int^{s_1}_0\psi^*_1(\tau)d\tau>\int^{s_1}_0 z^*_1(\tau)d\tau,\\
z^*_1(s), & s\in[0,s_1]\ \mbox{if}\ \int^{s_1}_0\psi^*_1(\tau)d\tau<\int^{s_1}_0 z^*_1(\tau)d\tau,\\
\psi^*_1(s), & s\in(s_1 ,s_0),\\
z^*_1(s), & s\in[s_0,|B_R(0)|].
\end{array}
\right.
\end{equation*}
Then, we can verify that $w(s)$ satisfies
$$
-\frac{dw(s)}{ds}\leq
\lambda_1(B_R,\frac{R^*}{R}\beta)n^{-2}w_n^{-\frac{2}{n}}s^{\frac{2}{n}-2}\int^s_0w(\tau)d\tau\ \ \ \mbox{for any}\ \ \ s\in(0,\ |B_R(0)|).
$$
At this stage, a similar argument to that of the case (i) can lead to a contradiction.

Summing up, we complete the proof of Theorem 3.1.

\section*{4.\ The proof of Theorem 1.1 and 1.6}

\setcounter{section}{4}

\setcounter{equation}{0}

This section devotes to prove Theorem 1.1 and Theorem 1.6. Some Lemmas needed in the proof of Theorem 1.1 are presented in the Appendix.
\vskip0.1cm
{\bf Proof of Theorem 1.1:}\  Let $z_1$ be the first eigenfunction, and $z_2$ be the radial part of the second eigenfunction, of problem (\ref{eq1.7}). Set
\begin{equation}\label{eq1.27.1}
g(r)=\left\{\begin{array}{ll}
\frac{z_2(r)}{z_1(r)}, & 0\leq r<R,\\
\lim\limits_{x\rightarrow R^{-}}g(r), & r\geq R,
\end{array}
\right.
\end{equation}
and
\begin{equation}\label{eq1.27.3}
\eta(r)={g'(r)}^2+\frac{n-1}{r^2}g^2(r).
\end{equation}
By Lemma A.2, we have
\begin{equation}
\lambda_2(\Omega,\beta)-\lambda_1(\Omega,\beta)\leq
\frac{\int_{\Om}\eta( r)\psi^2_1dx}{\int_{\Om}g^2( r)\psi^2_1dx}.
\label{eq1.18}
\end{equation}
Combining the conclusion of Lemma A.4 with (\ref{eq1.18}), we get
\begin{equation}
\lambda_2(\Omega,\beta)-\lambda_1(\Omega,\beta)\leq
\frac{\int_{B_R}\eta(r)z^2_1dx}{\int_{B_R}g^2(r)z^2_1dx}.
\label{eq1.30.2}
\end{equation}
It follows from Lemma A.5 and (\ref{eq1.30.2}) that
\begin{eqnarray}\label{eq4.1}
\lambda_2(\Omega,\beta)-\lambda_1(\Omega,\beta)\leq\lambda_2(B_R, \frac{R^*\beta}{R})-\lambda_1(B_R, \frac{R^*\beta}{R}).
\end{eqnarray}
By the rescaling property of eigenvalue problems, (\ref{eq4.1}) can be rewritten as
\begin{eqnarray}\label{eq4.2}
\lambda_2(\Omega,\beta)-\lambda_1(\Omega,\beta)\leq\frac{R_\lambda^2}{R^2}[\lambda_2(B_{R_\lambda}, \frac{R^*\beta}{{R_\lambda}}) -\lambda_1(B_{R_\lambda} \frac{R^*\beta}{{R_\lambda}})].
\end{eqnarray}
Again, by the rescaling property of eigenvalue problems, we have
\begin{eqnarray}\label{eq4.3}
\lambda_1(\Omega,\beta)=\lambda_1(B_{R_\lambda},\frac{R^*}{R_\lambda}\beta)
\end{eqnarray}
and
\begin{eqnarray}\label{eq4.4}
\lambda_i(B_{R_\lambda},\frac{R^*}{R_\lambda}\beta)=\frac{{R^*}^2}{R_\lambda^2}\lambda_i(\Om^*,\be)\ \ \mbox{for}\ \ i=1,2.
\end{eqnarray}
From (\ref{eq4.2}), (\ref{eq4.3}) and (\ref{eq4.4}), we can finally deduce that
$$\frac{\la_2(\Om,\be)}{\la_1(\Om,\be)}\leq\frac{R_\lambda^2}{R^2}\frac{\la_2(\Om^*,\be)}{\la_1(\Om^*,\be)}-\frac{R_\lambda^2}{R^2}+1.$$
This completes the proof of Theorem 1.1.

\vskip 0.2cm
{\bf Proof of Theorem 1.6:}\ Let $\psi_1(x)$ be the first eigenfunction of the problem (\ref{eq1.1}), and $z_1(x)$ be the first eigenfunction of the problem (\ref{eq1.7}). For any $p>0$, we set
$$f(x)=\psi_1(x)/||\psi_1||_{L^p(\Omega)},\ \ \mbox{and}\ \ g(x)=z_1(x)/||z_1||_{L^p(\Omega)}.$$
Obviously, $f(x)$ and $g(x)$ are also the first eigenfunction of the problem (\ref{eq1.1}) and the (\ref{eq1.7}) respectively. Moreover,
$$\int_{\Om}|f(x)|^pdx=\int_{B_R(0)}|g(x)|^pdx=1.$$
At this stage, we divide the proof of Theorem 1.6 into the following two cases.
\vskip 0.1cm
(i)\ In the case $g^*(|B_R|)<f^*(|B_R|)$, it follows from Theorm 3.1 that there exists a unique point $s_0\in(0,\ M)$ such that
\begin{equation*}
\left\{\begin{array}{ll}
f^*(s)\leq g^*(s), & s\in [0,\ s_0],\\
f^*(s)>g^*(s), & s\in ( s_0,\ |B_R(0)|].
\end{array}
\right.
\end{equation*}
From this, we can deduce that
$$\int^{s}_0|f^*(t)|^pdt\leq\int^{s}_0|g^*(t)|^pdt\ \ \ \mbox{for any}\ \ \ s\in[0,\ |B_R(0)|]$$
due to $\int^{|\Om|}_0|f^*(s)|^pds=\int^{|B_R(0)|}_0|g^*(s)|^pds=1$.

Since $f^*(s)$ and $g^*(s)$ are non-increasing, by applying Proposition 2.4 and Proposition 2.8 to $f^*(s)$ and $g^*(s)$, we conclude that
$$ \int^{|\Om|}_0|f^*(s)|^qds\leq\int^{|B_R(0)|}_0|g^*(s)|^qds\ \ \ \mbox{for any}\ \ \ q\geq p>0.$$
Hence
$$\int_{\Om}|f(x)|^qdx\leq\int_{B_R(0)}|g(x)|^qdx\ \ \ \mbox{for any}\ \ \ q\geq p>0.$$
By the definition of $f(x)$ and $g(x)$, we get
$$\left(\int_{\Om}|\psi_1|^qdx\right)^{\frac{1}{q}}\leq K(p,q,\beta,\Omega,n)\left(\int_{\Om}|\psi_1|^pdx\right)^{\frac{1}{p}}$$
with
 \begin{eqnarray*}
 K(p,q,\beta,\Omega,n)=\left(\int_{B_R(0)}|z_1|^qdx\right)^{\frac{1}{q}}\bigg/
 \left(\int_{B_R(0)}|z_1|^pdx\right)^{\frac{1}{p}}.
\end{eqnarray*}
This is just the desired conclusion of Theorem 1.6.
\vskip 0.1cm
(ii)\ In the case $g^*(|B_R|)\geq f^*(|B_R|)$, it follows from Theorem 3.1 that $g^*(s)\geq f^*(s)$ for
any $s\in[0,\ |B_R|]$. Thus, we have
$$\int^s_{0}|f^*(t)|^pdt\leq\int^s_{0}|g^*(t)|^pdt\ \ \ \mbox{for any}\ \ \ s\in[0,\ |B_R(0)|].$$
With this inequality, we can obtain the conclusion of Theorem 1.6 in a similar way to that of the case (i).

\section*{5.\ The Proof of Remark 1.5}

\setcounter{section}{5}

\setcounter{equation}{0}

\noindent
In this section, we give a sketch proof of Remark 1.5. To this end, we let
$$\Phi_1(x)=|\nabla \psi_1|^2+\frac{2}{n}\la_1(\Om,\be)\psi^2_1.$$
By a similar argument to that used in \cite{payne}, we can conclude that $\Phi_1(x)$ takes its maximum either on $\partial\Omega$, or at an interior point $p$ with $\nabla\psi_1(p)=0$. Moreover, by similar computations to that used in \cite{payne}, we can prove that $\Phi_1(x)$ can not take its maximum on $\partial\Om$ if $\Omega$ is convex. Hence, for any $\tilde{p}\in\partial\Omega$, it holds
$$
|\nabla \psi_1(\tilde{p})|^2+\frac{2}{n}\la_1(\Om,\be)\psi^2_1(\tilde{p})\leq\frac{2}{n}\la_1(\Om,\be)M^2_{\Om},
$$
with $M_{\Om}=\max\limits_{x\in\Om}\psi_1(x)$. Taking the boundary condition into account, we can get
$$
\left(\be^2+\frac{2}{n}\la_1(\Om,\be)\right)\psi^2_1(\tilde{p})\leq\frac{2}{n}\la_1(\Om,\be)M^2_{\Om}.
$$
Thus
\begin{equation}
M\leq\sqrt{\frac{\frac{2}{n}\la_1(\Om,\be)}{\be^2+\frac{2}{n}\la_1(\Om,\be)}}M_{\Om}.\label{5.1}
\end{equation}
From the proof of Lemma 3.2, we know that
$$
n^{2}w_n^{\frac{2}{n}}\leq-\lambda_1(\Omega,\beta)M_{\Om}\mu_{\psi_1}^{\frac{2}{n}-1}{\mu'_{\psi_1}}(t)\ \ \ \mbox{for any}\ \ \ t>M.
$$
Integrating the above inequality on $(M, M_{\Om})$, we have
\begin{eqnarray}\label{eq5.2}
n^{2}w_n^{\frac{2}{n}}(M_{\Om}-M)\leq\frac{n}{2}\lambda_1(\Omega,\beta)M_{\Om}|\Om_M|^{\frac{2}{n}}
\end{eqnarray}
due to $\mu_{\psi_1}(M_{\Om})=0$ and $\mu_{\psi_1}(M)=|\Om_{M}|$.

From (\ref{5.1}) and (\ref{eq5.2}), we have
\begin{eqnarray*}
|\Om_M|^{\frac{2}{n}}\geq\frac{2nw_n^{\frac{2}{n}}}{\lambda_1(\Omega,\beta)}
\left(1-\sqrt{\frac{\frac{2}{n}\la_1(\Om,\be)}{\be^2+\frac{2}{n}\la_1(\Om,\be)}}\right).
\end{eqnarray*}
Since $f(x)=\frac{x}{x+\be^2}$ is increasing on $(0,+\infty)$ and $\lambda_1(\Omega,\beta)\leq\lambda_1(\Omega)$, we get
$$R_M\geq\bigg[\frac{2n}{\lambda_1(\Omega)}\bigg(1-\sqrt{\frac{\frac{2}{n}\la_1(\Om)}
{\be^2+\frac{2}{n}\la_1(\Om)}}\bigg)\bigg]^{\frac{1}{2}}.$$

\vskip 0.1cm

\section*{6.\ Appendix}

\setcounter{section}{6}

\setcounter{equation}{0}

\noindent

In this appendix, we outline the proof of these lemmas used in the proof of Theorem 1.1 in section 4.
\vskip 0.1cm
{\em {\bf Lemma A.1.}\ For any $P(x)$ such that $P(x)\not\equiv 0$ and
$\int_{\Om}P\psi^2_1dx\equiv0$, we
$$
\lambda_2(\Omega,\beta)-\lambda_1(\Omega,\beta)\leq
\frac{\int_{\Om}|\nabla
P|^2\psi^2_1dx}{\int_{\Om}P^2\psi^2_1dx}.
$$

}
\vskip 0.1cm
{\bf Proof:}\ From the Rayleigh-Ritz inequality for $\la_2(\Omega,\beta)$, we have

$$\lambda_2(\Omega,\beta)\leq\frac{\int_{\Om}|\nabla u|^2dx
+\be\int_{\partial\Om}u^2dS}{\int_{\Om}u^2dx}
$$
for any $u$ satisfying $u\not\equiv 0$ and $\int_{\Om}u\psi_1dx=0$.

Taking $u=P\psi_1$ as a trial function, we obtain
$$
\lambda_2(\Omega,\beta)\leq\frac{\int_{\Om}|\nabla
P|^2\psi^2_1dx+2\int_{\Om}\psi_1P\nabla P\cdot \nabla
\psi_1+P^2|\nabla\psi_1|^2dx+\be\int_{\partial\Om}P^2\psi_1^2dS}{\int_{\Om}P^2\psi^2_1dx}.
$$
It is easy to check that
$$
\lambda_1(\Omega,\beta)=\frac{2 \int_{\Om}\psi_1P\nabla P\cdot \nabla \psi_1+P^2|\nabla\psi_1|^2dx+\be\int_{\partial\Om}P^2 \psi_1^2dS}{\int_{\Om}P^2\psi^2_1dx}.
$$
Hence
\begin{equation}
\lambda_2(\Omega,\beta)-\lambda_1(\Omega,\beta)\leq \frac{\int_{\Om}|\nabla P|^2\psi^2_1dx}{\int_{\Om}P^2\psi^2_1dx}.\label{eq1.12}
\end{equation}
\vskip 0.1cm
Let $g_0:\ (0,\ +\infty)\mapsto R$ be a nonnegative nontrivial bounded continuous differential function.
We consider the mapping $T:\ R^n\mapsto R^n$ which is defined by
\begin{equation}\label{eq1.13}
T(x_0)=\int_{\Om}g_0(|x-x_0|)\frac{x-x_0}{|x-x_0|}\psi^2_1dx.
\end{equation}
If $B$ is a ball containing $\Om$, then it is obvious that $T(x_0)$ points inward on $\partial B$. Hence, it follows from the Brouwer fixed point theorem that there exists $x^*_0\in B$ such that
$$
Tx^*_0=x^*_0.
$$
 Choosing $x^*_0$ as the origin of $R^n$ , we have
 \begin{equation}
\int_{\Om}g_0(r)\frac{x}{r}\psi^2_1dx=0.\label{eq1.14}
\end{equation}
Hence, $P_i=g_0(r)\frac{x_i}{r},\ i=1,2,\dots,n$ can be used as trial functions in Lemma A.1, and we can get a lemma as the following.
\vskip 0.1cm
 {\em {\bf Lemma A.2.}\ For any  nonnegative bounded continuous and differentiable function $g_0(r)$, it holds
 \begin{equation}
\lambda_2(\Omega,\beta)-\lambda_1(\Omega,\beta)\leq \frac{\int_{\Om}
[g'_0(r)^2+\frac{n-1}{r^2}g_0^2(r)]\psi^2_1dx}{\int_{\Om}g^2_0(r)\psi^2_1dx}.
\end{equation}

}

\vskip 0.1cm
{\bf Proof:}\ From (\ref{eq1.12}), we have
$$
[\lambda_2(\Omega,\beta)-\lambda_1(\Omega,\beta)]\int_{\Om}P^2\psi^2_1dx
\leq\int_{\Om}|\nabla P|^2\psi^2_1dx.
$$
Taking $P(x)=P_i=g_0(r)\frac{x_i}{r},\ i=1,2,\dots,n$, it yields
$$
[\lambda_2(\Omega,\beta)-\lambda_1(\Omega,\beta)]\int_{\Om}P_i^2\psi^2_1dx
\leq\int_{\Om}|\nabla P_i|^2\psi^2_1dx,\ \ i=1,2,\dots,n.
$$
Summing on $i$ for $i=1,2,\dots,n$, we obtain
\begin{equation}\label{eq6.1}
\lambda_2(\Omega,\beta)-\lambda_1(\Omega,\beta)\leq\frac{\int_{\Om}\sum^n_i|\nabla P_i|^2\psi^2_1dx}{\int_{\Om}\sum^n_iP^2_i\psi^2_1dx}.
\end{equation}
Since
$$
\sum^n_iP^2_i=g^2_0(r)\ \ \ \mbox{and}\ \ \ \sum^n_i|\nabla P_i|^2=(g'_0)^2+\frac{n-1}{r^2}g_0^2(r),
$$
it follows from (\ref{eq6.1}) that
\begin{equation}
\lambda_2(\Omega,\beta)-\lambda_1(\Omega,\beta)\leq \frac{\int_{\Om}
[g'_0(r)^2+\frac{n-1}{r^2}g_0^2(r)]\psi^2_1dx}{\int_{\Om}g^2_0(r)\psi^2_1dx}.
\end{equation}
This completes the proof of Lemma A.2.
\vskip 0.1cm

{\em {\bf Lemma A.3.}\ Let $g(r)$ and $\eta(r)$ be the functions given by (\ref{eq1.27.1}) and (\ref{eq1.27.3}) respectively. Then,
$g(r)$ is increasing and $\eta(r)$ is decreasing.}

{\bf Proof:}\ As in \cite{Ben}, we define a function $q(r):= \frac{rg'(r)}{g(r)}$ for any $r\in[0, R]$. Then, it is clear that the conclusion of Lemma A.3 is equivalent to $0\leq q\leq1$ and  $q'(r)\leq0$  for $0\leq r\leq R$. By the definition of $g(r)$, we can rewrite $q(r)$ as
\begin{equation}
q(r)=r(\frac{z'_2}{z_2}-\frac{z'_1}{z_1}).
\label{eq1.21}
\end{equation}
Since $z_1$ is the first eigenfunction, and $z_2$ is the radial part of the second eigenfunction, of problem (\ref{eq1.7}), it follows that $z_1(r)$ and $z_2(r)$ satisfy the following differential equations respectively.
\begin{equation}
z''_1+\frac{n-1}{r}z'_1+\lambda_1(B_R, R^*\beta/R)z_1=0
\label{eq1.19}
\end{equation}
\vskip 0.1cm
\begin{equation}
z''_2+\frac{n-1}{r}z'_2+\left(\lambda_2(B_R, R^*\beta/R)-\frac{n-1}{r^2}\right)z_2=0.
\label{eq1.20}
\end{equation}
Combining (\ref{eq1.21}), (\ref{eq1.19}) and (\ref{eq1.20}), we can show that, for any $r\in(0,R)$, $q(r)$ satisfies the following Riccati equation
\begin{equation}\label{eq6.2}
q'(r)=(\lambda_1(B_R, R^*\beta/R)-\lambda_2(B_R, R^*\beta/R))r+(1-q)(q+n-1)/r-2q\frac{z'_1}{z_1}.
\end{equation}

Let $J_p(x)$ denote the Bessel function of order $p$. Then, it is well known that $z_1(r)=Cr^{1-\frac{n}{2}}J_{\frac{n}{2}-1}(\sqrt{\lambda_1(\Omega, \beta)}r)$ for $0\leq r\leq R$. Therefore, by the property of Bessel function, the Riccati equation (\ref{eq6.2}) can be rewritten as
\begin{equation}\label{eq6.3}
\begin{array}{ll}
q'(r)&=(\lambda_1(B_R, R^*\beta/R)-\lambda_2(B_R, R^*\beta/R))r+(1-q)(q+n-1)/r\\
      &+2\sqrt{\lambda_1(\Omega, \beta)}q\frac{J_{\frac{n}{2}}(\sqrt{\lambda_1(\Omega, \beta)}r)}{J_{\frac{n}{2}-1}(\sqrt{\lambda_1(\Omega, \beta)}r)}.
\end{array}
\end{equation}

Before proceeding on, we first consider the behavior of $q(r)$ at the endpoints $r=0$ and $r=R$. Since $z_1(r)$  and $z_2(r)$ satisfy the following boundary conditions
$$ z_1(0)<+\infty,\ \ \ z'_1(0)=0,\ \ \ \frac{dz_1(R)}{dr}+\frac{R^*\beta}{R}z_1(R)=0$$
and
$$z_2(0)=0,\ \ \   \frac{dz_2(R)}{dr}+\frac{R^*\beta}{R}z_2(R)=0,$$
we can show, by L'Hoptital's rule, that
$$q(0)=1,\ \ \  q'(0)=0$$
and
$$q(R)=0,\ \ \ q'(R)=\left(\lambda_1(B_R, R^*\beta/R)-\lambda_2(B_R, R^*\beta/R)\right)R+(n-1)/R.$$

Now, we are in a position to prove $0\leq q(r)\leq1$ and $q'(r)\leq 0$ for $r\in[0,R]$. At first, we prove $q(r)\geq 0$ for $0\leq r\leq R$ by contradiction. To this end, we suppose in contrary that $q(r)$ changes sign in $[0,R]$. Then, from the facts that $q(0)=1$ and $q(R)=0$, we may conclude that there should exist two points $r_1$ and $r_2$ with $0< r_1<r_2\leq R$ such that $q(r_1)=q(r_2)=0$, $q'(r_1)\leq 0$ and $q'(r_2)\geq 0$. On the other hand, by the Riccati equation (\ref{eq6.2}), we have
\begin{eqnarray*}
0\geq q'(r_1)&=&(\lambda_1(B_R, R^*\beta/R)-\lambda_2(B_R, R^*\beta/R))r_1+(n-1)/r_1\\
             &>&(\lambda_1(B_R, R^*\beta/R)-\lambda_2(B_R, R^*\beta/R))r_2+(n-1)/r_2=q'(r_2)\geq 0
\end{eqnarray*}
which is a contradiction. Therefore, we have $q(r)\geq0$ for $0\leq r\leq R$.

In the second, we prove $q(r)\leq 1$ for $0\leq r\leq R$. Suppose in contrary. Then, there exists two points with $0< r_1<r_2<R$ such that $q(r_1)=q(r_2)>1$, $q'(r_1)\geq0$  and  $q'(r_2)\leq 0$. Since $J_{p+1}(x)/xJ_p(x)$ is strictly increasing on $[0,\sqrt{\lambda_1(B_1)})$ for $ p\geq-1/2$ (see\cite{Ben}), it follows from (\ref{eq6.3}) that
\begin{eqnarray*}
0\leq\frac{1}{r_1}q'(r_1)&=&(\lambda_1(B_R, R^*\beta/R)-\lambda_2(B_R, R^*\beta/R))+(1-q)(q+n-1)/r^2_1\\
&+&2\sqrt{\lambda_1(\Omega, \beta)}q\frac{J_{\frac{n}{2}}(\sqrt{\lambda_1(\Omega, \beta)}r_1)}{r_1J_{\frac{n}{2}-1}(\sqrt{\lambda_1(\Omega, \beta)}r_1)}\\
             &<&(\lambda_1(B_R, R^*\beta/R)-\lambda_2(B_R, R^*\beta/R))+(1-q)(q+n-1)/r^2_2\\
&+&2\sqrt{\lambda_1(\Omega, \beta)}q\frac{J_{\frac{n}{2}}(\sqrt{\lambda_1(\Omega, \beta)}r_2)}{r_2J_{\frac{n}{2}-1}(\sqrt{\lambda_1(\Omega, \beta)}r_2)}\\
             &=&\frac{1}{r_2}q'(r_2)\leq0,
\end{eqnarray*}
which is a contradiction. Thus, $q(r)\leq 1$ for any $r\in[0,R]$.

At last, we prove $q'(r)\leq 0$ for $0\leq r\leq R$.  Suppose not. We can find three points $r_1,$ $r_2,$ $r_3$ with $0< r_1<r_2<r_3<R$ such that $q(r_1)=q(r_2)=q(r_3)$, $q'(r_1)\leq 0$, $q'(r_2)\geq 0$, and $q'(r_3)\leq 0$. Writing $r_2$ as $r_2=tr_1+(1-t)r_3$ for some $t\in(0,1)$, and using the convexity of $\frac{1}{r}$, $(\lambda_1(B_R, R^*\beta/R)\\-\lambda_2(B_R, R^*\beta/R))r$, and $J_\frac{n}{2}(r)/J_{\frac{n}{2}-1}(r)$ (see\cite{Ben}), we obtain from the equation (\ref{eq6.3}) that
\begin{eqnarray*}
0\leq q'(r_2)&=&(\lambda_1(B_R, R^*\beta/R)-\lambda_2(B_R, R^*\beta/R))r_2+(1-q)(q+n-1)/r_2\\
                      &+& 2\sqrt{\lambda_1(\Omega, \beta)}q\frac{J_{\frac{n}{2}}(\sqrt{\lambda_1(\Omega, \beta)r_2})}{J_{\frac{n}{2}-1}(\sqrt{\lambda_1(\Omega, \beta)}r_2)}\\
                      &<&t q'(r_1)+(1-t)q'(r_3)\leq 0.
\end{eqnarray*}
This is a contradiction. Hence, $q'(r)\leq 0$ for any $r\in[0,R]$, and the proof of Lemma A.3 is completed.

 \vskip 0.1cm
 {\em {\bf Lemma A.4.}\ Let $g(r)$ and $\eta(r)$ be the functions given by (\ref{eq1.27.1}) and (\ref{eq1.27.3}) respectively. If we normalize $\psi_1(x)$ and $z_1(x)$ so that $\int_{\Om}\psi^2_1(x)dx=\int_{Om}z^2_1(x)dx$, then there holds
\begin{equation}
 \int_{\Om}\eta(r)\psi^2_1dx\leq \int_{B_R(0)}\eta( r)z^2_1dx\label{eq1.22}
\end{equation}
 and
\begin{equation}
 \int_{\Om}g^2( r)\psi^2_1dx\geq\int_{B_R(0)}g^2( r)z^2_1dx.\label{eq1.23}
\end{equation}
}
 {\bf Proof:}\ By Lemma A.3 and the properties of rearrangement, we have
  \begin{equation}
 \int_{\Om}\eta( r)\psi^2_1dx\leq \int_{\Om^*}\eta( r){\psi^\star_1}^2dx\label{eq1.24}
\end{equation}
and
 \begin{equation}
 \int_{\Om}g^2( r)\psi^2_1dx\geq \int_{\Om^*}g^2( r){\psi^\star_1}^2dx.\label{eq1.25}
\end{equation}
Hence, in order to prove Lemma A.4, we only need to prove
 \begin{equation}
 \int_{\Om^*}\eta( r){\psi^\star_1}^2dx\leq\int_{B_R(0)}\eta( r)z^2_1dx\label{eq1.26}
\end{equation}
and
\begin{equation}
\int_{\Om^*}g^2( r){\psi^\star_1}^2dx\geq \int_{B_R(0)}g^2( r)z^2_1dx.\label{eq1.27}
\end{equation}

In the case $z^*_1(|B_R(0)|)<\psi^*_1(|B_R(0)|)$, if we set $r_1=\left(\frac{s_0}{w_n}\right)^{\frac{1}{n}}$, then it follows from Theorem 3.1 that
\begin{eqnarray*}
&&\int_{B_R(0)}\eta( r)z^2_1dx-\int_{\Om^*}\eta(r){\psi^{\star}_1}^2dx\\
&=&nw_n\bigg[\int^{r_1}_0\eta( r)(z^2_1-{\psi^{\star}_1}^2)r^{n-1}dr\\
& &+\int^{R}_{r_1}\eta(r)(z^2_1-{\psi^{\star}_1}^2)r^{n-1}dr-\int^{R^*}_R\eta(r){\psi^{\star}_1}^2r^{n-1}dr\bigg]\\
&\geq& nw_n\eta( r_1)\left[\int^{r_1}_0(z^2_1-{\psi^{\star}_1}^2)r^{n-1}dr+
\int^{R}_{r_1}(z^2_1-{\psi^{\star}_1}^2)r^{n-1}dr-
 \int^{R^*}_R{\psi^{\star}_1}^2r^{n-1}dr\right]\\
 &=&\eta(r_1)\left[\int_{B_R(0)}z^2_1dx-
 \int_{\Om^*}{\psi^{\star}_1}^2dx\right]\\
 &=&0.
\end{eqnarray*}
and
\begin{eqnarray*}
&&\int_{B_R(0)}g( r)z^2_1dx-\int_{\Om^*}g^2( r){\psi^{\star}_1}^2dx\\
                                   &=&nw_n\bigg[\int^{r_1}_0g(r)(z^2_1-{\psi^{\star}_1}^2)r^{n-1}dr\\
                                   & &+
                                   \int^{R}_{r_1}g( r)(z^2_1-{\psi^{\star}_1}^2)r^{n-1}dr-
                                   \int^{R^*}_R g( r){\psi^{\star}_1}^2r^{n-1}dr\bigg]\\
                                   &\leq&nw_n g(r_1)\left[\int^{r_1}_0(z^2_1-{\psi^{\star}_1}^2)r^{n-1}dr+
                                   \int^{R}_{r_1}(z^2_1-{\psi^{\star}_1}^2)r^{n-1}dr-
                                   \int^{R^*}_R{\psi^{\star}_1}^2r^{n-1}dr\right]\\
                                   &=&g(r_1)[\int_{B_R(0)}z^2_1dx-\int_{B_R(0)}{\psi^{\star}_1}^2dx]\\
                                   &=&0.
\end{eqnarray*}
This is just the conclusion we want.
\vskip 0.1cm
In the case $z^*_1(|B_R(0)|)\geq\psi^*_1(|B_R(0)|)$, we can prove the conclusion of Lemma A.4 in a similar way as above by making use of Theorem 3.1.
\vskip 0.1cm

 {\em {\bf Lemma A.5.}\ Let $g(r)$ and $\eta(r)$ be the functions given by (\ref{eq1.27.1}) and (\ref{eq1.27.3}) respectively. If $z_1$ is the first eigenfunction, and $z_2$ is the radial part of the second eigenfunction, of problem (\ref{eq1.7}), then we have
\begin{equation}\label{eq1.28}
\lambda_2(B_R, \frac{R^*\beta}{R})-\lambda_1(B_R,
\frac{R^*\beta}{R})=\frac{ \int_{B_R}\eta(r)z^2_1dx}{
\int_{B_R}g^2(r)z^2_1dx}.
\end{equation}
}

{\bf Proof:}\ Multiplying the equation (\ref{eq1.19}) by $\frac{z_2^2}{z_1}$ and then integrating  on $B_R(0)$, we have
\begin{eqnarray}\label{eq1.31}
\lambda_1(B_R,\frac{R^*\beta}{R})\int_{B_R(0)}z_2^2
=\frac{R^*\beta}{R}\int_{\partial B_R(0)}z_2^2+\int_{B_R(0)}
\left(\frac {z_2^2}{z_1}\right)^{'} {z_1}^{'}
\end{eqnarray}
by virtue of the boundary conditions $z_1(0)<+\infty$, $z'_1(0)=0$ and $\frac{\partial
z_1(R)}{\partial\nu}+\frac{R^*\beta}{R} z_1(R)=0$.

Multiplying the equation (\ref{eq1.20}) by $z_2$ and then
integrating on $B_R(0)$, we have
\begin{eqnarray}\label{eq1.32}
\lambda_2(B_R,\frac{R^*\beta}{R})\int_{B_R(0)}z_2^2
=\int_{B_R(0)}( {z_2}^{'})^2+\frac{R^*\beta}{R}\int_{\partial
B_R(0)}z_2^2+\int_{B_R(0)}\frac{n-1}{r^2}z_2^2
\end{eqnarray}
by virtue of the boundary condition $z_2(0)=0$ and $\frac{\partial
z_2(R)}{\partial\nu}+\frac{R^*\beta}{R} z_2(R)=0 $.

Obviously, the conclusion of Lemma A.5 can be deduced from (\ref{eq1.31}) and (\ref{eq1.32}) immediately.


\vskip 0.5cm
\begin{thebibliography}{s2}

\bibitem{Ben}M.S.Ashbaugh, R.D. Benguria, proof of the Payne-P$\acute{o}$lya-Weinberger conjecture, Bull of A.M.S, 25(1991), 19-29.


\bibitem{Ashb}M.S.Ashbaugh, R.D. Benguria, A sharp bound for the ratio of the first two eigenvalues of Dirichlet Laplacians and extensions, Ann. Math., 135(1992), 601-628.

\bibitem{Ashbaugh eigenvalue ratios}M. S. Ashbaugh and R. D. Benguria, More bounds on eigenvalue ratios for
Dirichlet Laplacians in n dimensions, SIAM J. Math. Anal., 24 (1993), 1622-
1651.

\bibitem{Ashbaugh-Neumann}M.S.Ashbaugh, R.D.Benguria ,Universal Bounds for the Low Eigenvalues of Neumann Laplacians in N Dimensions, SIAM J. Math. Anal., 24(3),1993 ,557-570.


\bibitem{Bandle}C.Bandle, Isoperimetric inequalities and applications, Pitman Publishers,London(1980).


\bibitem{ref GC}G. Chiti, A reverse H$\ddot{o}$lder inequality
                for the eigenfunctions of linear second order
                elliptic operators, Journal of Applied Mathematics and
                Physics, 33 (1982), 143-148.

\bibitem{G. Chiti2}G. Chiti, An isoperimetric inequality for the eigenfunctions of linear second order elliptic
operators. Boll. Un. Mat. Ital. A (6) 1, 145-151 (1982).

\bibitem{Bos} M.H.Bossel, Membranes $\acute{e}$lastiquement
li$\acute{e}$es inhomog$\grave{e}$nes ou sur une surface: une
nouvelle extension du th$\acute{e}$or$\acute{e}$me
isop$\acute{e}$rim$\acute{e}$trique de Rayleigh-Faber-Krahn, Z.
Angew. Math. Phys., 39(1988), 733-742.

\bibitem{sIssa} Chavel Isaac, Eigenvalues in Riemannian Geometry,
Academic Press, Inc, Orlando, FL, 1984.


\bibitem{QD}Qiuyi Dai and Yuxia Fu, Faber-Krahn inequality for Robin problem involving p-
                Laplacian, Acta Mathematicae Applicatae Sinica, English Series, 27(2011), 13-28.


\bibitem{Dane1}D.Daners, A Faber-Krahn inequality for the Robin problems in any space dimension, Math. Ann, 335(2006), 767-785.

\bibitem{Dane2} D.Daners, J. Kennedy, Uniqueness in the Faber-Krahn inequality
for Robin problems, SIAM J. Math. Anal., 39(2007), 1191-1207.


\bibitem{Fab} G.Faber, Beweis, dass unter allen homogenen Membranen
von gleicher Fl\"{a}che und gleicher Spannung die kreisf\"{o}rmige
den tiefsten Grundton gibt, Sitzungsber, Bayr. Akad. Wiss.
M\"{u}nchen, Math.-Phys. Kl, 1923, 169-127.

\bibitem{Gio} T.Giorgi, R.G.Smits, Monotonicity results for the
principal eigenvalue of the generalized robin problem, Illinois J.
Math. 49 (4) (2005), 1133-1143.


\bibitem{ref GJG} G. H. Hardy, J. E. Littlewood and G. P$\acute{o}$lya,
Some simple inequalities satisfied by convex functions. Messenger
Math.58(1929), p152.


\bibitem{Henrot1}A.Henrot, Minimization problems for eigenvalues of the Laplacian, J.Eovl Equation,3(2003)443-461.

\bibitem{Henrot2}A.Henrot, E.Oudet, Minimizing the second eigenvalue of the Laplace operator with Dirichlet boundary conditions, Arch. Rational Mech.Anal, 169(2003),73-87.

\bibitem{ref BK}B. Kawohl, Rearrangements and convexity of Level sets in
PDES, Lecture notes in Mathematics, 1150, Springer-Verlag,
Heidelberg 1985.

\bibitem{ref SK}S. Kesavan, Symmetrization and Applicantions. Series in Anylysis, Vol.3,
           World Scientific Books, April,2006.

\bibitem{sKeP} S.Kesavan, F.Pacella, Symmetry of positive
solutions of a quasilinear elliptic equation via isoperimetric
inequalities, Applicable Anal., 54(1994), 27-37.

\bibitem{Kra} E.Krahn, \"{U}ber eine von Rayleigh formulierte
Minimaleigenschaft des Kreises, Math. Ann. 94 (1925), 97-100.

\bibitem{Kroger} P.Kr$\ddot{o}$ger, Upper bounds for the Neumann eigenvalues on a bounded domains in Euclidean space, J.Funct. Anal.106(1992),353-357.

\bibitem{PPW1}L.E. Payne, G.P$\acute{o}$lya,H.F.Weinberger, Sur le quotient de devx fr$\acute{e}$quences propres cons$\acute{e}$cutives, Comptes Rendus Acad. Sci. Paris, 241(1955)917-919.

\bibitem{PPW2}L.E. Payne, G.P$\acute{o}$lya,H.F.Weinberger, On the ratio of consecutive eigenvalues, J.Math. and Phys.35(1956).

\bibitem{PayRay}L.E. Payne, M.E.Rayner, An isoperimetric inequality for the first eigenfunction in the fixed membrane problem, Z. Angew. Math. Phys, 23(1972)13-15.

\bibitem{payne}L.E. Payne, P.W.Schaefer, Eigenvalue and eigenfunction inequalities for the elastically supported membrane, Z.Angew.Math.Phys, 52(2001)888-895.

\bibitem{PS}L.E. Payne, I.Stakgold, On the mean value of the fundamental mode in the fixed membrane problem. Applicable Analysis, 3(1973)295-303.


\bibitem{szego}G.Szeg$\ddot{o}$, Inequalities for certain eigenvalues of a membrane of given area, J.Ratioal Mech.Anal, 3(1954), 343-356.

\bibitem{Takahashi}F.Takahashi and A.Uegaki, A Payne-Rayner Type Inequality for the
Robin Problem on Arbitrary Minimal Surfaces in $R^n$, Results in Mathematics,
59(2011), 107-114.



\bibitem{Jobin}M.Th$\acute{e}$re$\grave{s}$e and K.Jobin, isoperimetric monotonicity and Isoperimetric inequalities of Payne-Rayner type
for the first eigenfunction of the Helmholtz problem, J. Appl.Math. Phys.,32 (1981),625-646.


\bibitem{Thompson}C.J.Thompson, on the ratio of conseccutive eigenvalues in N-dimensions, Stud.Appl.Math.48(1969),281-283.

\bibitem{QWang}Qiaoling Wang and Changyu Xia, Isoperimetric bounds for the first eigenvalue of
the Laplacian, Z. Angew. Math. Phys. 61 (2010), 171-175.

\bibitem{Weinberger}H.F.Weinberger, An isoperimetric inequality for the $n-$dimensional free membrane problem, J.Ratioal Mech.Anal, 5(1956), 633-636.















\end{thebibliography}
\end{document}